\documentclass[twoside, 11pt]{article}
\usepackage{quoting}
\usepackage[letterpaper,margin=1in]{geometry}
\usepackage{enumitem}
\usepackage{amsthm}
\usepackage{tikz}
\usetikzlibrary{decorations.pathreplacing,calligraphy}

\usepackage{amsmath}

\usepackage[linktocpage=true,pagebackref=true,colorlinks,linkcolor=magenta,citecolor=blue,bookmarks,bookmarksopen,bookmarksnumbered,ocgcolorlinks]{hyperref}
\usepackage{comment}
\usepackage{lineno}
\usepackage{array}
\usepackage[nameinlink]{cleveref}
\usepackage[nottoc]{tocbibind}
\usepackage{pifont}
\usetikzlibrary{calc}
\usepackage{hyperref}
\usepackage{amsthm, amsmath, amssymb}
\usepackage{mathrsfs}
\usepackage{mathtools}
\usepackage{ragged2e}
\usepackage{algpseudocode,algorithm}
\usepackage{cleveref}
\usepackage{lineno}
\usepackage{comment}
\usepackage{todonotes}
\usepackage{amsthm}
\usepackage{tikz}
\usepackage{circuitikz}
\usetikzlibrary{decorations.pathreplacing,calligraphy}
\usepackage{caption}
\usepackage{subcaption}
\usepackage{wrapfig}
\usepackage{amsmath}

\usepackage[charter]{mathdesign}

\newtheorem{definition}{Definition}
\newtheorem{observation}{Observation}
\newtheorem{claim}{Claim}

\newtheorem{theorem}{Theorem}
\newtheorem{lemma}{Lemma}
\newtheorem{corollary}{Corollary}

\newcommand{\T}{TCC(a_1,a_2,\dots,a_{s-1})}
\newcommand{\TT}{TCC(a_1,a_2,\dots,a_s)}
\newcommand{\M}{TCC(m_1,m_2,\dots,m_d)}
\newcommand{\TCH}[1]{TC(H_ {#1})}
\newcommand{\TTCH}[1]{TC^*(H_ {#1})}
\newcommand{\Y}{\Upsilon}

\newcommand{\G}{\Gamma}

\newcommand{\Z}{\mathbb{Z}}

\newcommand{\ceil}[1]{\lceil #1 \rceil}

\providecommand{\keywords}[1]
{
  \small	
  \textbf{\textit{Keywords---}} #1
}

\def\centerarc[#1](#2)(#3:#4:#5)
    { \draw[#1] ($(#2)+({#5*cos(#3)},{#5*sin(#3)})$) arc (#3:#4:#5); }

\tikzset{
    cross/.pic = {
    \draw[rotate = 45] (-#1,0) -- (#1,0);
    \draw[rotate = 45] (0,-#1) -- (0, #1);
    }
}

\usepackage{titling}
\predate{}
\postdate{}

\title{Boxicity and Cubicity of A Subclass of Divisor Graphs and Power Graphs of Cyclic Groups}
\date{}
 \author{
L. Sunil Chandran \\
 Indian Institute of Science, India \\
    \texttt{sunil@iisc.ac.in} \\
  \and
    Jinia Ghosh \\
  IIT Gandhinagar, India \\
  \texttt{jiniag@iitgn.ac.in} \\  
}
\everypar{\looseness=-1}

\begin{document}

\maketitle




\begin{abstract}
The \textit{boxicity} (\textit{cubicity}) of an undirected graph $\Gamma$ is the smallest non-negative integer $k$ such that $\Gamma$ can be represented as the intersection graph of axis-parallel rectangular boxes (unit cubes) in $\mathbb{R}^k$. An undirected graph is classified as a \textit{comparability graph} if it is isomorphic to the comparability graph of some partial order. This paper studies boxicity and cubicity for subclasses of comparability graphs.

We initiate the study of boxicity and cubicity of a special class of algebraically defined comparability graphs, namely the \textit{power graphs}. The power graph of a group is an undirected graph whose vertex set is the group itself, with two elements being adjacent if one is a power of the other. We analyse the case when the underlying groups of power graphs are cyclic. Another important family of comparability graphs is \textit{divisor graphs}, which arises from a number-theoretically defined poset, namely the \textit{divisibility poset}. We consider a subclass of divisor graphs, denoted by $D(n)$, where the vertex set is the set of positive divisors of a natural number $n$, and two vertices $a$ and $b$ are adjacent if and only if $a$ divides $b$ or $b$ divides $a$.

We first show that to study the boxicity (cubicity) of the power graph of the cyclic group of order $n$, it is sufficient to study the boxicity (cubicity) of $D(n)$. We derive estimates, tight up to a factor of $2$, for the boxicity and cubicity of $D(n)$. The exact estimates hold good for power graphs of cyclic groups.
\end{abstract}

\keywords{Boxicity, Cubicity, Poset Dimension, Comparability Graphs, Divisor Graphs, Power Graphs}

\section{Introduction}\label{intro}

Let $\Gamma$ be a simple undirected graph. We write $V(\G)$ and $E(\G)$ to denote the \textit{vertex set} and the \textit{edge set} of $\G$, respectively. For basic definitions and notations from graph theory, an interested reader can refer to the books by West \cite{west2001introduction} or any other standard textbook. An edge of $E(\G)$ between $u$ and $v$ is denoted by $\{u,v\}$. A \textit{universal} vertex in a graph is a vertex which is adjacent to all other vertices of the graph. 
 An undirected graph $I$ is called \textit{interval graph} (\textit{unit interval graph}) if there is a mapping $f$ from $V(I)$ to the set of closed (unit) intervals on the real line $\mathbb{R}$ such that two vertices $u$ and $v$ are adjacent in $I$ if and only if the corresponding intervals intersect, i.e., $f(u) \cap f(v) \neq \emptyset$. The mapping $f$ is called an \textit{interval representation} of $I$.

\subsection{Boxicity and Cubicity}

\begin{definition}\label{box cub def}\textbf{(Boxicity, Cubicity)}
     A $k$-box ($k$-cube) is the Cartesian product of $k$ closed (unit) intervals on the real line $\mathbb{R}$. A $k$-box ($k$-cube) representation of a graph $\G$ is a function $f$ that maps a vertex of $\G$ to a $k$-box ($k$-cube) such that two distinct vertices $u$ and $v$ are adjacent in $\G$ if and only if $f(u) \cap f(v) \neq \emptyset$. The boxicity (cubicity) of $\G$, denoted by $box(\G)$ $(cub(\G))$, is the smallest non-negative integer $k$ such that $\G$ has a $k$-box ($k$-cube) representation.
     
\end{definition}

From the above definition, it is clear that the boxicity (cubicity) of complete graphs is $0$.
Unlike many other well-known graph parameters such as chromatic number, connectivity, or treewidth, boxicity and cubicity are not monotonic. That is, if $\G'$ is a supergraph of $\G$, boxicity (cubicity) of $\G'$ may be greater or less than that of $\G$, depending on $\G$ and $\G'$.  However,  for induced subgraphs, the following inequality holds (that follows from the definition): 

\begin{lemma}\label{subgraph}
 If $\G'$ is an induced subgraph of a graph $\G$, then $box(\G') \leq box(\G)$ and $cub(\G') \leq cub(\G)$.
 \end{lemma}

A graph has boxicity (cubicity) at most 1 if and only if it is an (unit) interval graph.  Note that a $k$-cube representation of $\G$ using unit cubes is equivalent to a $k$-cube representation where all cubes have a constant side length. Clearly, for any finite graph $\Gamma$, $box(\G) \leq cub(\G)$.

The concepts of boxicity and cubicity were introduced by Roberts \cite{MR0252268} in 1969. 
He showed that every graph with $n$ vertices has a $\lfloor n/2 \rfloor$-box   ($\lfloor 2n / 3 \rfloor$-cube) representation. Computing the boxicity of a graph was shown to be NP-hard by Cozzens \cite{cozzens1982higher}. Later, Yannakakis \cite{yannakakis1982complexity} and Kratochvil \cite{kratochvil1994special} strengthened it, and the latter showed that deciding
whether the boxicity of a graph is at most $2$ is NP-complete. Yannakakis \cite{yannakakis1982complexity} proved that deciding whether the cubicity of a given graph is at least $3$ is NP-hard.

Several studies in the literature have explored the dependencies of boxicity and cubicity on the structural properties of graphs. Esperet \cite{esperet2009boxicity} proved that $box(\G) \leq \Delta^2 + 2$, where $\Delta$ denotes the maximum degree of the graph $\G$. Scott and Wood \cite{scott2020better} improved this result to $O(\Delta \log^{1+o(1)} \Delta)$. 
On a different note, Chandran et al. \cite{chandran2007boxicity} established that the boxicity of a graph is at most its tree-width plus 2. Chandran et al. \cite{chandran2009cubicity} derived bounds on both the boxicity and cubicity of a graph in terms of the size of the minimum vertex cover $t$, showing that $cub(\G) \leq t + \lceil \log(n - t) \rceil - 1$, and $box(\G) \leq \lceil t/2 \rceil + 1$.

Additionally, boxicity and cubicity have been studied for special graph classes, leading to the establishment of bounds for these parameters. Roberts \cite{MR0252268} proved that the boxicity of a complete $k$-partite graph is $k$.  Scheinerman \cite{scheinerman1984intersection} showed that the outerplanar graphs have boxicity at most $2$, while Thomassen \cite{thomassen1986interval} established that the planar graphs have boxicity bounded above by $3$. Felsner and Mathew \cite{felsner2011contact} provided a different proof for Thomassen's result using contact representations with cubes. Chandran et al. \cite{chandran2007boxicity} provided upper bounds for the boxicity of several special graph classes, including chordal graphs, circular arc graphs, AT-free graphs, permutation graphs, and co-comparability graphs. Similarly, the cubicity of certain graph classes, such as hypercubes and complete multipartite graphs, has been analyzed by Roberts \cite{MR0252268}, Chandran et al. \cite{chandran2005cubicity, chandran2008cubicity}, and Michael et al. \cite{michael2006sphericity}. Chandran, Mathew and Rajendraprasad \cite{chandran2016upper} established that $cub(\G) \leq
2 \ceil{\log_2 \chi (\G)} box(\G) + \chi (\G) \ceil{\log_2 \alpha(\G)}$, where $\chi(\G)$ and $\alpha(\G)$ denote the chromatic number and independence number of $\G$, respectively. Recently, Cauduro et al. \cite{caoduro2023boxicity} provided bounds for the boxicity of Kneser graphs and of the complement of the line graph of any graph. Boxicity has a connection with another parameter called separation dimension \cite{ alon2015separation,basavaraju2016separation} as Basavaraju et al. \cite{basavaraju2016separation} showed that the separation dimension of a hypergraph is equal to the boxicity of its line graph. Majumder and Mathew \cite{majumder2022local} provided upper bounds for the local boxicity of a graph in terms of its maximum degree, the number of edges and the number of vertices.

Algorithmic approaches for computing boxicity have also been explored. Adiga et al. \cite{adiga2010parameterized} initiated an algorithmic study and provided parameterized algorithms for boxicity using parameters such as vertex cover and max-leaf number. Furthermore, Adiga, Babu and Chandran et al. \cite{adiga2018sublinear} developed polynomial time $o(n)$-factor approximation algorithms for boxicity and related parameters, including poset dimension and interval dimension.

\subsection{Partial Order and Poset Dimension}

A \textit{partially ordered set or poset} is an ordered pair $P=(S,\leq_P)$, where $S$ is a set and $\leq_P$ is a \textit{partial order}, i.e, $\leq_P$ is a reflexive, anti-symmetric, transitive binary relation defined on the set $S$, known as the \textit{ground set}. For $a,b \in S$, we write $a < _P b$
to mean that $a \leq_P b$ and $a \neq b$. Two elements $a$ and $b$ are said to be \textit{comparable} if $a \leq_P b$ or $b\leq_P a$, otherwise $a$ and $b$ are called as \textit{incomparable}. A \textit{linear order} is a poset where the elements are pairwise comparable. 

A \textit{linear extension} $L=(S,\leq_L)$ of a partial order $P=(S,\leq_P)$ is a linear order which satisfies the following: if $x \leq_P y$, then $x \leq_L y$. A \textit{realizer} $\mathcal{R}$ of a poset $P = (S, \leq_P)$ is a set of linear
extensions of $P$ that follows the condition: for any two distinct elements $x$ and $y$ in $S$, $x <_P y$ if and only if $x <_L y$  for all $L \in R$ (This means that if $x$ and $y$ are incomparable in $P$, then there exist $L_1,L_2 \in  \mathcal {R} $ such that $x <_{L_1} y$  and $y <_{L_2} x$.). The \textit{poset dimension}
of $P$, denoted by $dim(P)$, is the minimum positive
integer $k$ such that there exists a realizer of $P$ of cardinality $k$. Poset dimension, introduced by Dushnik and Miller (1941) \cite{dushnik1941partially}, is a fundamental and extensively studied concept of partial orders (see \cite{trotter1992combinatorics, trotter1995partially} for more references).

The \textit{comparability graph} $\G_P$ \textit{of a poset} $P=(S,\leq_P)$ is an undirected graph with $S$ as the vertex set and $\{x,y\} \in E(\G_P)$ if and only if $x$ and $y$ are comparable in $P$. Moreover, an undirected graph $\G$ is a comparability graph if there is a poset $P$ such that its corresponding comparability graph $\G_P$ is isomorphic to $\G$. Note that it is possible for the same comparability graph to be isomorphic to the comparability graphs of more than one distinct partial order.  But as shown by Trotter et al. \cite{trotter1976dimension}, all 
the posets that give rise to the same comparability graph, have the same poset dimension, throwing light on the fundamental nature of the parameter, poset dimension.

 In 1982, when Yannakakis proved the NP-completeness of both poset dimension and boxicity \cite{yannakakis1982complexity}, he might have intuitively suspected a potential connection between the boxicity of a comparability graph and the poset dimension of its underlying poset.
 In 2011, a link between these two concepts is established as Adiga et al. \cite{adiga2011boxicity} proved the following result:

\begin{lemma}[\cite{adiga2011boxicity}]\label{poset dim and boxicity}
Let $\G_P$ be the comparability graph of a poset $P$. Then the following holds: $\frac{box(\G_P)}{\chi(\G_P)-1} \leq dim(P) \leq 2\cdot box(\G_P)$, where $\chi(\G_P)$ denotes the chromatic number of $\G_P$.
\end{lemma}

Adiga et al. \cite{adiga2011boxicity} discussed the intimate relation of partial order dimension and boxicity in their paper. From their discussion, it seems that despite its seemingly geometric definition, boxicity is `almost' a combinatorial parameter, undeniably with a tinge of geometry in it. Thus, the study of boxicity for subclasses of comparability graphs is as natural as the study of partial order dimension for the corresponding partial orders. On the other hand, cubicity seems to be closer to geometry; for example, allowing us to use the total volume of the cubes associated with the set of vertices in an independent set of graph $\G$ to derive lower bounds for cubicity of $\G$.

\subsection{Power Graphs}\label{sec: power graphs}

A group is a set with a binary operation which satisfies certain constraints. The basic definitions and facts from group theory can be found in any standard book (e.g.,  \cite{rotman2012introduction}). In this paper, we only consider \emph{finite groups}.

\begin{definition}\label{def power graph}\textbf{(Power graph)}
    The \emph{power graph} of a group $G$, denoted by $Pow(G)$, is an undirected graph with vertex set $G$, where two distinct vertices $x$ and $y$ are adjacent if $x=y^m$ or $y=x^m$ for some positive integer $m$.
\end{definition}

The notion of power graphs was originally defined by Chakrabarty et al. \cite{chakrabarty2009undirected}. For surveys on power graphs,  interested readers can refer to \cite{abawajy2013power} and \cite{kumar2021recent}. It is known that the power graphs are comparability graphs \cite{feng2014power, aalipour2017structure}. Note that, if $b = a^k$ and $c=b^\ell$, then $c = a^{kl}$; therefore, it is possible to transitively orient the edges of any power graph. This makes the poset dimension and the boxicity (cubicity) interesting parameters to examine in power graphs.

\subsection{Divisor Graph and the Transitive Closure of Cartesian Product of Complete Graphs}\label{sec: divisor graphs}

A divisor graph is the comparability graph of a number theoretically defined poset, namely divisibility poset. A \textit{divisibility poset} is a poset $P_D=(S,\leq)$, where $S$ is a subset of positive integers and the order relation is defined as follows: for any two elements $a,b \in S$, $a \leq b$ if and only if $a \mid b$. The comparability graph of a divisibility poset is referred to as a \textit{divisor graph}. Though the name  `divisor graphs' has started gaining popularity only recently  \cite{santhosh1828divisor}, divisibility posets were among the most well-known partially ordered sets, appearing as even textbook material right from the inception of the study of posets. Moreover, it is easy to show that any finite poset is isomorphic to some divisibility poset. To see this, consider a poset $P=(S,\leq_P)$ with ground set $S=\{x_1,x_2,\dots,x_n\}$ and $n$ distinct primes $p_1,p_2,\dots, p_n$. With the mapping $f(x_i)={ \prod_{i_k \in I_k} p_{i_k}}$, where $I_k$ is the set of indices $i_j$ for which $x_{i_j} \leq x_i$, an isomorphic divisibility poset with ground set $f(S) = \{ f(x) : x \in S\}$  can be constructed from $P$. Therefore, the class of divisor graphs is the same as the class of comparability graphs. Thus, the study of the poset dimension of divisor graphs is the study of the poset dimension in its most general setting. Hence, it makes sense to
consider sub-classes of divisor graphs with more structure. 

In this paper, we study a special class of finite divisor graphs defined as follows: for a natural number $n$, let $D(n)$ denote the divisor graph where the vertex set is the set of divisors of $n$.  In \Cref{divisor and power}, we explain how the boxicity of this sub-class of divisor graphs is related to the boxicity of power graphs \textit{of groups of a special class - class of finite cyclic groups.} 

We now digress a little and define \textit{the transitive closure of the Cartesian product of complete graphs}, which we will use to investigate the boxicity and cubicity of the class of divisor graphs $D(n)$.

\vspace{0.3cm}

\noindent \textbf{Transitive closure of the Cartesian product of complete graphs:} The \textit{Cartesian product} of two
graphs $\G_1$ and $\G_2$, denoted  by $\G_1 \square \ \G_2$, is a graph on the vertex set
$V(\G_1) \times V(\G_2)$ with the following edge set: $E(\G_1 \ \ \square \ \G_2) = \{\{(u_1, u_2),(v_1, v_2)\} \ : \ u_1 = v_1, \{u_2,v_2\} \in E(\G_2) \text{ or }
\{u_1,v_1\} \in E(\G_1), u_2 = v_2\}$. We define the Cartesian product of more than two graphs as follows: For $s > 2$,  $\square_{i=1}^{s} \G_i = (\square_{i=1}^{s-1} \G_i) \square \G_s$.  Let the vertex set of a complete graph $K_q$ be $\{0,1,2,\dots,q-1\}$, where $q \geq 2$. For positive integers $a_1,a_2, \ldots, a_s$, define  $ C(a_1,a_2, \ldots, a_s) = \square_{i=1}^{s}K_{a_i+1}$. Clearly, $C(a_1,\ldots,a_s)$  has the vertex set  $\{(x_1,x_2,\dots,x_s) \ : \ x_i \in \Z \text{ and } 0 \leq x_i \leq a_i, \text{ for all } 1 \leq i \leq s\}$. For an $s$-tuple $x$, we write $x_i$ to denote the $i$th component of $x$. Thus,  $x=(x_1,x_2,\dots,x_s)$. Two distinct vertices $x$ and $y$ (which are also $s$-tuples) are adjacent in $C(a_1,a_2,\dots,a_s)$ if and only if there is exactly one $k$ in $[s]$ such that $x_k \ne y_k$  and for all $i \in [s]\setminus \{k\}$, $x_i=y_i$.

Given that the vertices of each constituent complete graph \( K_q \) are labeled as \( \{0, 1, 2, \dots, q-1\} \), we can define a binary relation \( \prec \) on \( C(a_1, a_2, \dots, a_s) \) as follows: for any two distinct elements \( x \) and \( y \) in \( C(a_1, a_2, \dots, a_s) \), we say \( x \prec y \) if and only if there exists exactly one \( k \in [s] \) such that \( x_k < y_k \) (where \( < \) is the standard ordering on \( \mathbb{Z} \)) and \( x_i = y_i \) for all \( i \in [s] \setminus \{k\} \). Observe that if \( \{x, y\} \) is an edge in \( C(a_1, a_2, \dots, a_s) \), then either \( x \prec y \) or \( x \succ y \). Thus, the binary relation \( \prec \) provides an orientation for the edges of \( C(a_1, a_2, \dots, a_s) \), where the edge \( \{x, y\} \) is assigned the direction \( (x, y) \) if \( x \prec y \). The graph obtained by taking the transitive closure of this oriented graph of $C(a_1,a_2,\dots,a_s)$ is denoted by \( TCC(a_1, a_2, \dots, a_s) \) and is referred to as the \textit{transitive closure of the Cartesian product of complete graphs} \( K_{a_1+1}, K_{a_2+1}, \dots, K_{a_s+1} \). For convenience, the corresponding comparability graph (the underlying undirected graph) of $\TT$ is also called by the same name and denoted by the same notation. In the remaining discussion of the paper, we use $\TT$ to refer to its corresponding comparability graph.


Note that the vertex set of the graph $TCC(a_1,a_2\dots,a_s)$ is the set $\{x=(x_1,x_2,\dots,x_s): \ x_i \in \Z \text{ and } 0 \leq x_i \leq a_i, \text{for all } i \in [s]\}$ and two distinct vertices ($s$-tuples) are adjacent in $TCC(a_1,a_2,\dots,a_s)$ if either $x_i \leq y_i$, for all $1\leq i \leq s$ (referred to as $x \leq y$) or $x_i \geq y_i$, for all $1\leq i \leq s$ (referred to as $x\geq y$). 
If two distinct $s$-tuples $x$ and $y$ in $TCC(a_1,a_2,\dots,a_s)$ satisfy one of these conditions, they are referred to as comparable; otherwise, they are considered incomparable. 

\noindent {\bf Special cases of Cartesian products of complete graphs:} 
\begin{enumerate}

\item {\bf The case when $a_1 =a_2 = \ldots = a_s$:} The \textit{Hamming graph} is the Cartesian product of complete graphs of the same size, i.e., in our definition, if $a_i=q-1$, for all $i \in [s]$ we get a Hamming graph which is the Cartesian product of $s$-copies of complete graphs of $q$ vertices. 

\item  {\bf The case when $a_1 =a_2 = \ldots = a_s = 1$:} If we consider the Cartesian product of $s$ copies of $K_2$ (complete graphs of $2$ vertices), we get an \textit{$s$-dimensional hypercube} $H_s$, a very popular graph in graph theory and combinatorics. Let $TC(H_s)$ denote the \textit{transitive closure of $H_s$}.  It is well-known that  $TC(H_s)$ is exactly the comparability graph of the partial order defined on the power set of $[s]$, where the underlying relation is set inclusion. An $s$-tuple $x$ with component values from $\{0,1\}$ can be mapped to a subset $S_x$ of $[s]$ by the following rule: for each $i \in [s]$, 
$x_i = 1$ if and only if $i \in S_x$. The study of the poset dimension of this poset (which is also known as the Boolean lattice) has been a subject of extensive study since 1950; see \cite{dushnik1950concerning}, \cite{spencer1972minimal}, \cite{brightwell1994dimension} and \cite{kostochka1997dimension}. Therefore, this special case is of particular interest.  Understanding the boxicity of this comparability graph $TC(H_s)$ is crucial in our proofs.


Note that $\mathbf{0}=(0,0,\dots,0)$ and $\mathbf{1}=(1,1,\dots,1)$ are the universal vertices in $\TCH{s}$. Let $\TTCH{s}$ denote the graph $\TCH{s}\setminus \{\mathbf{0},\mathbf{1}\}$ and we call it the \textit{truncated transitive closure of hypercube}. It is easy to see that the removal of universal vertices does not change the boxicity of a graph since universal vertices can be mapped to very large boxes of dimension $d$ that contain all other boxes, where $d$ is the boxicity of the remaining graph obtained after removing the universal vertices. 

\begin{observation}
$box(\TCH{s})=box(\TTCH{s})$.
\end{observation}

\end{enumerate}

We make an easy observation that the divisor graph $D(n)$ is isomorphic to the transitive closure of the Cartesian product of certain complete graphs. This, in turn, shows that the structure of  $D(n)$ does not depend on the values of the primes in the prime factorization of $n$; it depends only on the number of primes involved and their exponents.

\begin{observation}\label{equivalence of divisor graph and tcc}
    If $n=p_1^{a_1}p_2^{a_2}\dots p_s^{a_s}$, then the divisor graph $D(n)$ is isomorphic to $TCC(a_1,a_2,
    \dots,a_s)$. 
\end{observation}
\begin{proof}
    We define a mapping $f$ from the vertex set of $D(n)$ to the vertex set of $TCC(a_1,a_2,\dots,a_s)$ as follows:
    $f(p_1^{x_1}p_2^{x_2}\dots p_s^{x_s})=(x_1,x_2,\dots,x_s)$. Here $0 \leq x_t \leq a_t$ and for each $t \in [s]$. It is easy to verify that $f$ is indeed an isomorphism between $D(n)$ and $TCC(a_1,a_2,
    \dots,a_s)$ since for two distinct numbers  $n_1 = p_1^{x_1}p_2^{x_2}\dots p_s^{x_s}$ and $n_2 = p_1^{y_1}p_2^{y_2}\dots p_s^{y_s}$, $n_1 \mid n_2$ if and only if $(x_1,x_2,\dots,x_s) < (y_1,y_2,\dots,y_s)$.
\end{proof}

Recall that a number is squarefree if all the exponents in its prime factorisation are 1. Based on the above observation, it is clear that for a squarefree number $n$ with $d$ distinct prime factors, the divisor graph $D(n)$ is isomorphic to the transitive closure of the $d$-dimensional hypercube $TC(H_d)$. Consequently, determining the boxicity (cubicity) of the divisor graph $D(n)$ for a squarefree number $n = p_1p_2 \dots p_d$ is equivalent to determining the boxicity (cubicity) of $TC(H_d)$. On a similar note,  this observation suggests that finding the boxicity (cubicity) of $D(n)$, for $n$ whose prime factors are raised to the same power, is equivalent to 
finding the boxicity (cubicity) of the transitive closure of Hamming graphs. These special cases are noteworthy for two reasons: (1)  The poset dimension of the underlying poset associated with $ TC(H_d)$ was extensively studied; therefore, it is interesting to study its boxicity also.  (2)  There are some papers by Chandran et al. \cite{chandran2008cubicity} and Imrich et al. \cite{chandran2015boxicity}, which establish bounds on the cubicity and boxicity of hypercubes and Hamming graphs without considering their transitive closures. From this context, it may be interesting to explore the boxicity (cubicity) of the transitive closures of the above-mentioned graphs.



\subsection{Divisor Graph $D(n)$ and Power Graph of Cyclic Group of Order $n$}\label{divisor and power}

  A subset $H$ of a group $G$ is called a \textit{subgroup} of $G$ if $H$ forms a group under the binary operation of $G$. This is denoted by $H \leq G$. The number of elements in a group $G$ is called the \textit{order of the group}, denoted by $|G|$. \textit{Lagrange's theorem} says that for every subgroup $H$ of a finite group $G$, $|H|$ divides $|G|$. 

The \textit{order of an element} $g$ in $G$, denoted by $o(g)$, is the smallest positive integer $m$ such that $g^m=e$, where $e$ is the identity element of $G$. For an element $g \in G$ with $o(g)=m$, the set $\{g, g^2, \dots, g^{m-1}, g^m=e\}$ forms a subgroup in $G$; it is called the \textit{cyclic subgroup generated by} $g$ and denoted by $\langle g \rangle$.
A group $G$ is called \emph{cyclic} if $G=\langle g \rangle$ for some $g\in G$. Such an element $g$ is called a \textit{generator} of $G$. A finite cyclic group of order $n$ is isomorphic to $\mathbb{Z}_n$, the group of integers additive modulo $n$. Moreover, a finite cyclic group of order $n$ has a unique subgroup (which is also cyclic) of order $d$ for each positive divisor $d$ of $n$.

\vspace{0.2cm}

\noindent \textit{Reduced Power graph}: 
An equivalence relation $\sim$ can be defined on a group $G$ as follows: for $x,y\in G$, $x\sim y$ if and only if $\langle x \rangle = \langle y \rangle$, i.e., $x$ and $y$ generate the same cyclic subgroup in $G$.  The equivalence class containing $x$ under $\sim$ is denoted by $[x]$. Note that $[x]$ is the set of generators of the cyclic subgroup $\langle x \rangle$. So, the elements of an equivalence class $[x]$ generate one another and hence $[x]$ makes a clique in $Pow(G)$. Moreover, all the elements of an equivalence class are of the same order. We define the \emph{order of an equivalence class} by the order of any element belonging to the class. 

For two distinct equivalence classes $[x]$ and $[y]$, if an element $x\in [x]$ is adjacent to an element $y\in [y]$ in $Pow(G)$, then every element of $[x]$ is adjacent to every element of $[y]$, since each element in $[x]$ and $[y]$ generate one another.  For a finite group $G$, we define an undirected graph $R(G)$ where the vertex set is the set of equivalence classes of $G$, and we add/introduce an edge between two distinct vertices $[x]$ and $[y]$, if and only if $\{x,y\}$ is an edge in $Pow(G)$. We call $R(G)$ the reduced power graph of $G$. 

\begin{observation}\label{boxicity of reduced power graph}
    $box(R(G))=box(Pow(G))$ and $cub(R(G))=cub(Pow(G))$.
\end{observation}
\begin{proof}
    It is easy to note that $R(G)$ is isomorphic to an induced subgraph of $Pow(G)$. Hence, by \Cref{subgraph}, 
    $box(R(G)) \leq box(Pow(G))$. For the opposite direction, let $f$ be a $k$-box representation of $R(G)$. We define a box representation $h$ of $Pow(G)$ as follows: for an element $x$ in $G$, $h(x)= f([x])$. We will prove that $h$ is indeed a $k$-box representation of $Pow(G)$. Note that each equivalence class $[x]$ induces a clique in $Pow(G)$ and for any two elements $x_1,x_2 \in [x]$, $h(x_1) \cap h(x_2) = f([x])\neq \emptyset$. By the construction, two distinct equivalence classes (vertices in $R(G)$) $[x]$ and $[y]$ are adjacent in $R(G)$ if and only if there is an edge between $x$ and $y$ in $Pow(G)$. Hence $x$ and $y$ are adjacent in $Pow(G)$ if and only if $h(x)\cap h(y)=f([x]) \cap f([y]) \neq \emptyset$. 

    Similar proof holds for cubicity.
\end{proof}

We have already seen that the power graph is a comparability graph (see \Cref{sec: power graphs}), and each comparability graph is isomorphic to a divisor graph (see \Cref{sec: divisor graphs}). Hence, a power graph is a divisor graph. But the converse is not true. However, if we restrict the vertex set of 
the divisor graph to be the set of divisors of a natural number $n$, the resulting divisor graph $D(n)$ is isomorphic to the reduced power graph of the cyclic group of order $n$ (noting that the cyclic group of order $n$ unique up to isomorphism). It is likely that researchers in the field of power graphs are already familiar with this connection. However, to the best of our knowledge, no formal proof exists in the literature. Therefore, we provide a proof here to make it accessible to all readers and for completeness.

\begin{observation}\label{reduced power graph is a divisor graph}
 The reduced power graph $R(G)$ of the finite cyclic group $G$ of order $n$ is isomorphic to the divisor graph $D(n)$ of $n$.
\end{observation}
\begin{proof}
    Since the finite cyclic group $G$ of order $n$ has a unique subgroup (which is also cyclic) of order $d$ for each positive divisor $d$ of $n$, there is a bijective mapping, say $f$, from the vertex set of $R(G)$ to the set of positive divisors of $n$. Moreover, the mapping $f$ is defined as follows: for a vertex $[x]$ in $R(G)$, $f([x])=o(x)$, where $o(x)$ is the order of $x$ in $G=$ the order of the equivalence class $[x]=$   the order of the subgroup $\langle x \rangle$, which is the unique subgroup of order $o(x)$ in $G$.

    Now we prove that $f$ is an isomorphism between $R(G)$ and $D(n)$. Let $\{[x],[y]\}$ be  an edge in $R(G)$. We have to show that either $o(x) \mid o(y)$ or $o(y) \mid o(x)$. Since $\{[x],[y]\}$ is an edge, either $x$ generates $y$ or $y$ generates $x$, i.e., either $\langle x \rangle  \leq \langle y \rangle $ or $\langle y \rangle  \leq \langle x \rangle $. Hence, by Lagrange's theorem, either $o(x) \mid o(y)$ or $o(y) \mid o(x)$. For the opposite direction, assume that $\{f([x]),f([y])\}$ is an edge in $D(n)$. Hence, either $o(x) \mid o(y)$ or $o(y) \mid o(x)$. Suppose without loss of generality that $o(x) \mid o(y)$. Since $\langle y \rangle$ is a cyclic subgroup of order $o(y)$ and $o(x)$ is a divisor of $o(y)$,  the cyclic subgroup $\langle y \rangle $ contains a cyclic subgroup of order $o(x)$. Since the subgroup of order $o(x)$ is unique in $G$, namely $\langle x \rangle$,  $\langle x \rangle $ is a subgroup of $\langle y \rangle$. This means $y$ generates $x$ and $\{x,y\}$ is an edge in $Pow(G)$. Thus, by construction, $\{[x],[y]\}$ is an edge in $R(G)$.
\end{proof}


\subsection{Boxicity of Divisor Graph $D(n)$ from Well-known Tools}\label{intro: box of divisor graph}

 Recently, numerous researchers such as Lewis et al. \cite{lewis2021order}, Haiman \cite{haiman2023dimension}, and Souza et al. \cite{souza2024improved}, have conducted extensive studies on poset dimensions of divisibility posets for various ground sets. If the ground set is the set of divisors of a natural number $n$, then the poset dimension of the divisibility poset equals the number of distinct prime divisors of $n$ (see \cite{lewis2021order}).

\begin{lemma}\cite{lewis2021order}\label{dim of divisibility}
If $n=p_1^{a_1}p_2^{a_2}\dots p_s^{a_s}$, then the poset dimension of the divisibility poset $P_{D(n)}$ defined on the set of divisors of $n$ is $s$.
\end{lemma}

We provide an outline of the proof for the above lemma by constructing a realizer $\mathcal{R}$ of the divisibility poset $P_{D(n)}$. The realizer $\mathcal{R}$ consists of $s$ many linear extensions $L_1,L_2,\dots,L_s$, where each $L_i$ is constructed by concatenating linear orders in a specific pattern. To define these, we first consider the subsets \(S_{i,j}\) of the ground set \(S\) of \(P_{D(n)}\), where \(1 \leq i \leq s\), \(0 \leq j \leq \alpha_i\), and \(S_{i,j} = \{x \in S \mid \text{the power of } p_i \text{ in } x \text{ is } j\}\). Let $L_{i,j}$ denote a linear order on $S_{i,j}$, arranged by increasing values of the elements. Note that if $a \mid b$ and $a \neq b$, then $a$ is strictly less than $b$. Thus, each $L_{i,j}$ is a linear extension of the subposet of $P_{D(n)}$ restricted to $S_{i,j}$. Finally, for each $1 \leq i \leq s$, we define the linear order $L_i$ as follows: $L_{i,0} < L_{i,1} < L_{i,2}< \dots < L_{i,\alpha_i}$. It is easy to see that $\mathcal{R}=\{L_1,L_2,\dots,L_s\}$ is indeed a realizer of $P_{D(n)}$.

Note that, to derive an upper bound for boxicity of $D(n)$ using \Cref{poset dim and boxicity}, we need to know the chromatic number $\chi (D(n))$. In a comparability graph, the chromatic number is equal to the clique number (the number of vertices in the largest clique), and the clique number is the same as the number of vertices in the longest chain in the underlying poset. In $P_{D(n)}$, where $n=p_1^{a_1}p_2^{a_2}\dots p_s^{a_s}$, the number of vertices in the longest chain is $(a_1+a_2+\dots+a_s+1)$. Hence, using \Cref{poset dim and boxicity} and \Cref{dim of divisibility}, the following upper bound for boxicity can be derived:
 
 \begin{equation}\label{eq: upper bound from well-known}
     box(D(p_1^{a_1}p_2^{a_2}\dots p_s^{a_s})) \leq s \cdot (\sum_{i=1}^s a_i).
 \end{equation} 
 
 From the result by Chadran et al. \cite{chandran2009upper}, the known upper bound of cubicity in terms of boxicity is as follows: $cub(\G) \leq box(\G).\ceil{\log |V(\G)|}$. Thus, the following upper bound for the cubicity of $D(n)$ is achievable from \Cref{eq: upper bound from well-known}:
 
\begin{equation} \label{eq : upper bound for cubicity from well-known}
     cub(D(p_1^{a_1}p_2^{a_2}\dots p_s^{a_s})) \leq s \cdot (\sum_{i=1}^s a_i) \cdot \ceil{\log |V(D(n))|} 
 \end{equation}

 Even if we use the upper bound for cubicity provided by Adiga et al. \cite{adiga2010cubicity}, $cub(\G) \leq box(\G).\ceil{\log \alpha(\G)}$, where $\alpha(\G)$ is the independence number of $\G$, it can be shown that the above upper bound for $cub(D(n))$ (\Cref{eq : upper bound for cubicity from well-known}) will not improve significantly. 
 The upper bound for $cub(D(n))$ that we establish in this paper is better than the one mentioned in \Cref{eq : upper bound for cubicity from well-known}, by an approximate factor of 
 $\frac{1}{s. \ceil{\log |V(D(n))|}}$ (see \Cref{upper bound} for details). Moreover, since $box(\G) \leq cub(\G)$, our derived upper bound for $cub(D(n))$ also serves as an upper bound for $box(D(n))$, which is again an improvement from the easily achievable one (mentioned in \Cref{eq: upper bound from well-known}) by an approximate factor $\frac{1}{s}$.

The lower bound for boxicity of $D(n)$ that can be attained by combining \Cref{poset dim and boxicity} and \Cref{dim of divisibility}  is as follows:
 \begin{equation}\label{eq: low bound from well-known}
   box(D(p_1^{a_1}p_2^{a_2}\dots p_s^{a_s})) \geq \frac{s}{2}.  
 \end{equation}

 Our lower bound for the boxicity of $D(n)$ involves the powers of the prime factors of $n$ (see \Cref{lower bound} for details). This gives a better lower bound for $box(D(n))$ than the one mentioned in \Cref{eq: low bound from well-known}.

\subsection{Our Results}
In \Cref{equivalence of divisor graph and tcc}, we have already seen that $D(p_1^{m_1}p_2^{m_2}\dots p_d^{m_d})$ is isomorphic to $\M$, i.e., the transitive closure of the cartesian product of complete graphs. We mainly refer to the latter in the rest of our discussion. Note that, $TCC(m_1)$ is isomorphic to a complete graph and therefore $cub(TCC(m_1))= box(TCC(m_1)) = 0$. Recall that, $box(\G) \leq cub (\G)$ for any graph $\G$.

\begin{theorem}\textbf{(Upper Bound)}\label{upper bound}
    Let $m_1 \leq m_2 \leq \dots \leq m_d$ be natural numbers, where $d \geq 2$. Then the following holds:\\
   $box(\M) \le cub(\M) \leq m_1 + m_2 + \dots + m_{d-1}$.
\end{theorem}

Since $TC(H_d)$ is the same as $\M$ with $m_1=m_2=\dots=m_d=1$, the following is a straightforward corollary of \Cref{upper bound}. 

\begin{corollary}\label{hypercube}
 $box (TC(H_d)) \le cub(TC(H_d))\leq d-1$.
\end{corollary}





\begin{theorem}\textbf{(Lower Bound)}\label{lower bound}
Let $m_1 \leq m_2 \leq \dots \leq m_d$ be natural numbers, where  $d \geq 2$.

\vspace{0.2cm}

\noindent (1) If $d=2$, then $cub(TCC(m_1,m_2)) \geq box(TCC(m_1,m_2)) \geq m_1$.

\vspace{0.2cm}

\noindent(2) If $d\geq 3$, then the following holds:
\begin{equation} \label{eq1}
\textstyle
\begin{aligned}
    cub(\M) &\geq box(\M) \\
            &\geq (m_1 + m_3 + \dots + m_{d-2}) + m_{d-1}, \text{when } d \text{ is odd}, \\
            &\geq (m_2 + m_4 + \dots + m_{d-2}) + m_{d-1}, \text{when } d \text{ is even}.
\end{aligned}
\end{equation}

\end{theorem}

Note that the lower bound for the boxicity (cubicity) of $\M$ (\Cref{lower bound}) is at least half of the upper bound (\Cref{upper bound}) of the boxicity (cubicity). Furthermore, based on \Cref{upper bound} and \Cref{lower bound}, it can be observed that the upper and lower bounds for the boxicity (as well as cubicity) of $\M$ are the same when $d \leq 3$. Specifically, this can be summarized as follows:

\begin{corollary}\label{corollary: exact boxicity}
    Let $m_1 \leq m_2 \leq m_3$ be three natural numbers. Then, we have the following:\\
    (1) $box(TCC(m_1))=cub(TCC(m_1))=0$;\\
    (2) $box(TCC(m_1,m_2))=cub(TCC(m_1,m_2))=m_1$;\\
    (3) $box(TCC(m_1,m_2,m_3))=cub(TCC(m_1,m_2,m_3))=m_1+m_2$.
\end{corollary}

\section{Preliminaries}\label{preli}
Let $\G$ be a graph and $\G_i$, $1 \leq i \leq k$, be graphs on the same vertex set as $\G$ such that $E(\G) = \cap_{i=1}^{k}E(\G_i)$. Then, $\G$ is said to be the \textit{intersection} of graphs $\G_i$'s for $1 \leq i \leq k$ and denote it as $\G = \cap_{i=1}^{k}\G_i$. Recall the definitions of boxicity and cubicity (\Cref{box cub def}). We now present alternate definitions of boxicity and cubicity \cite{MR0252268}, which are due to the following lemma by Roberts \cite{MR0252268}.

\begin{lemma}\cite{MR0252268} The boxicity (cubicity) of a graph $\G$ is at most $k$, if there exist $k$ (unit) interval graphs $I_1, I_2, \dots, I_k$ such that $V(I_i)=V(\G)$ for each $i \in [k]$ and $\cap_{i=1}^k E(I_i) = E(\G)$.
    
\end{lemma}

Note that if $\G = \cap_{i=1}^k I_i $, then each $I_i$ is a supergraph of $\G$ and for every pair of vertices $u, v \in V(\G)$ with $\{u,v\} \notin E(\G)$, there exists some $i \in [k]$ such that $\{u,v\} \notin E(I_i)$. Hence, finding a $k$-box ($k$-cube) representation of a graph $\G$ is the same as finding $k$ (unit) interval supergraphs of $\G$ with the property that every pair of non-adjacent vertices in $\G$ is non-adjacent in at least one of those (unit) interval supergraphs.

\begin{definition}\label{boxicity def 2}\textbf{(Alternate definition of boxicity, cubicity)} The \emph{boxicity (cubicity)} of a non-complete graph $\G$ is the minimum positive integer $k$ such that $\G = \cap_{i=1}^k I_i $, where $I_i$'s are (unit) interval graphs.
    
\end{definition}

The \textit{join} of two graphs $\G_1$ and $\G_2$ on disjoint vertex sets, denoted by $\G_1 \vee \G_2$, is the graph
with vertex set $V(\G_1) \cup V(\G_2)$ and edge set $E(\G_1) \cup E(\G_2) \cup \{\{x_1, x_2\} : x_1 \in V(\G_1), x_2 \in V(\G_2)\}$. The next lemma shows that boxicity is additive under the join operation on graphs. The correctness of the lemma follows from \Cref{boxicity def 2} and the fact that an interval graph cannot contain an induced 4-cycle. For details, one can refer to \cite[Lemma 3]{trotter1979characterization}.  

 \begin{lemma}\label{join}
 $box(\G_1 \vee \G_2) = box(\G_1) + box(\G_2)$.
\end{lemma}

The join of $n$ vertex disjoint graphs $\G_1,\G_2,\dots, \G_n$, denoted by $\bigvee_{i=1}^n \G_i$, is the graph with the vertex set $\bigcup_{i=1}^n V(\G_i)$ and with the edge set $\big(\bigcup_{i=1}^n E(\G_i)\big) \cup \big(\bigcup_{1\leq i < j \leq n} \{ \{v_i,v_j\} : \ v_i \in V(\G_i), \  v_j \in V(\G_j)\}\big)$. With an inductive proof, one can generalise the statement of \Cref{join} as follows.

\begin{lemma}\label{general join}
    $box(\bigvee_{i=1}^n \G_i)=\sum_{i=1}^n box(\G_i)$.
\end{lemma}

\section{Proof of \Cref{upper bound}}

First, we prove the following lemma, which is the foundation for the main proof.

\begin{lemma}\label{inductive process}
    $cub(\TT) \leq cub(\T)+a_s$, where $s \geq 2$.
\end{lemma}
\begin{proof} Let \( cub(\T) = k \). The proof involves constructing \( k + a_s \) unit interval graphs \( I_1, I_2, \dots, I_k,\) \(I_{k+1}, I_{k+2}, \dots, I_{k+a_s} \) such that \(\TT\) is the intersection of these unit interval graphs. That means our aim here is to show that each unit interval graph $I_i$, $1 \leq I \leq k+a_s$, is a supergraph of $\TT$, and every pair of non-adjacent vertices in $\TT$ is non--adjacent in one of the unit interval graphs. The first \( k \) unit interval graphs \( I_1, I_2, \dots, I_k \) are referred to as \textit{Type 1 unit interval graphs}, while the remaining \( a_s \) unit interval graphs \( I_{k+1}, I_{k+2}, \dots, I_{k+a_s} \) are called \textit{Type 2 unit interval graphs}. For this proof, a pair of non-adjacent vertices in $\TT$ is referred to as a \textit{non-edge} in \(\TT\). The classification of the \( k + a_s \) unit interval graphs corresponds to two distinct types of non-edges in \(\TT\), which are explained in detail in the remainder of the proof.

Clearly $V(\TT)=\{ (x,0),(x,1),\dots,(x,a_s) \ : \ x \in V(\T)\}$. By definition of $E(\TT)$, it is easy to see that two distinct vertices $(x,b)$ and $(y,c)$ are adjacent in $\TT$ if and only if $x \leq y$ and $a \leq b$.

     Let $\{(x,b),(y,c)\}$ be a non-edge in $\TT$. Note that, in this case, $(x,b)$ and $(y,c)$ are two distinct vertices in $\TT$, and moreover $x \neq y$. We categorize the non-edges of $\TT$ into the following two types:

\vspace{0.2cm}

    \textit{Type 1 non-edge}: Let $\{x,y\}$ be a non-edge in $\T$. That means $x$ and $y$ are two incomparable $(s-1)$-tuples. Then, irrespective of the value of $b$ and $c$, $(x,b)$ and $(y,c)$ are incomparable $s$-tuples and hence $\{(x,b),(y,c)\}$ is a non-edge in $\TT$. Note that when $s=2$, there is no non-edge of Type 1.

    \textit{Type 2 non-edge:} Let $\{x,y\}$ be an edge in $\T$. That means $x$ and $y$ are two distinct comparable $(s-1)$-tuples. Without loss of generality, assume that $x < y$. Then, there exists at least an index $i$ in $1 \leq i \leq s-1$ such that $x_i < y_i$. Let $b > c$. Since the $s$-th indices of $(x,b)$ and $(y,c)$ satisfy $b > c$, $\{(x,b),(y,c)\}$ is a non-edge in $\TT$.

    \vspace{0.2cm}

    We now define the Type 1 unit interval graphs \( I_1, I_2, \dots, I_k \). Since $cub(\T)=k$, \Cref{boxicity def 2} implies that $\T$ is the intersection of $k$ unit interval graphs, say $I'_1,I'_2,\dots, I'_k$. For each $1\leq i \leq k$, let $h_i$ be an interval representation of the unit interval graph $I'_i$. For each $1\leq i \leq k$, we define an interval representation $f_i$ of $I_i$ as follows: 

   \begin{center}
       $f_i((x,b)) = h_i(x)$, for all $x \in V(\T)$ and for all $0 \leq b \leq a_s$.
   \end{center} 

    \begin{claim}\label{type 1 interval graphs are supergraphs}
        For each $I_i$, $1 \leq i \leq k$, $E(\TT) \subseteq E(I_i)$.
    \end{claim}
    \begin{proof}
        
        Let $(x,b)$ and $(y,c)$ be adjacent in $\TT$. Then, $(x,b)$ and $(y,c)$ are comparable $s$-tuples. Without loss of generality, we assume that $x \leq y$ and $b \leq c$. This implies that either $x=y$ or $\{x,y\} \in E(\T)$. Note that in the both cases, $h_i(x) \cap h_i(y) \neq \emptyset$, for all $1 \leq i \leq k$. Moreover, $f_i((x,b))=h_i(x)$ and $f_i((y,c))=h_i(y)$. Hence, $f_i((x,b)) \cap f_i((y,c)) \neq \emptyset$ and therefore $(x,b)$ and $(y,c)$ are adjacent in each $I_i$.
    \end{proof}

   \begin{claim}\label{type 1 interval graphs and non-edge}
       For any Type 1 non-edge $\{(x,b),(y,c)\}$ of $\TT$,
    there exists some Type 1 unit interval graph $I_l$, such that $\{(x,b),(y,c)\} \notin E(I_l)$.
   \end{claim} 
    \begin{proof}
        Since $\{(x,b),(y,c)\}$ is a non-edge of Type 1, $\{x,y\}$ is a non-edge in $\T$. Hence there exists an unit interval graph $I'_l$ in $I'_1,I'_2,\dots,I'_k$ such that $\{x,y\} \notin E(I'_l)$, i.e., $h_l(x)\cap h_l(y)=\emptyset$. So, by our construction, $f_l((x,b))\cap f_l((y,c))=\emptyset$ and hence $\{(x,b),(y,c)\} \notin E(I_l)$.
    \end{proof}

From \Cref{type 1 interval graphs are supergraphs} and \Cref{type 1 interval graphs and non-edge}, one can see that the intersection $\bigcap_{i=1}^k I_i$ of the Type 1 unit interval graphs contains each edge of $\TT$, but no Type 1 non-edge of $\TT$.  
However, each non-edge $\{(x,b),(y,c)\}$ of Type 2 is present in all Type 1 unit interval graphs of $\TT$. It is because in this case $\{x,y\} \in E(\T)$ and due to our construction, $f_i((x,b)) \cap f_i((y,c)) =h_i(x) \cap h_i(y) \neq \emptyset$, for all $1 \leq i \leq k$. To tackle this situation, we construct the Type 2 unit interval graphs $I_{k+1}, I_{k+2}, \dots, I_{k+a_s}$. Let $a_1+ a_2 + \dots + a_{s-1}=S$. For $x \in V(\T)$, let $|x|=\sum_{i=1}^{s-1}x_i$. Let $f_{k+a}$ be an interval representation of $I_{k+a}$, $1 \leq a \leq a_s$ defined as follows:

\begin{center}
    $f_{k+a}((x,a')) = [|x|,S+|x|]$, for all $0 \leq a' \leq a-1$,

     \hspace{1.4cm}        $= [|x|-S,|x|]$, for all $a \leq a' \leq a_s$.

\end{center}


     



Note that all the closed intervals in the interval representations defined above (for Type 2 unit interval graphs) have the same length $S$, and therefore, each closed interval can be scaled to a unit closed interval.



\begin{claim}\label{Type 2 interval graphs are supergraphs}
    For each $I_{k+a}$, $1 \leq a \leq a_s$, $E(\TT) \subseteq E(I_{k+a})$.
\end{claim}
\begin{proof}
    Let $x$ be a vertex in $V(\T)$. It is easy to note that $|x| \in \{0,1,2,\dots,S\}$. For each $(x,a')$, where $0\leq a' \leq a-1$, the point $S \in f_{k+a} ((x,a'))$. Hence, the set $\{(x,a') \ : \ x \in V(\T), \ 0 \leq a' \leq a-1\}$ induces a clique in $I_{k+a}$. Similarly, for each $(x,a')$, where $a\leq a' \leq a_s$, the point $0 \in f_{k+a} ((x,a'))$. Hence, the set $\{(x,a') \ : \ x \in V(\T), \ a \leq a' \leq a_s\}$ induces a clique in $I_{k+a}$.

    Now, we consider the only remaining case. Let $(x,b)$ and $(y,c)$ are adjacent in $\TT$, where $x \leq y$, $0 \leq b \leq a-1$ and $a \leq c \leq a_s$. Note that, $f_{k+a}((x,b))=[|x|,S + |x|]$ and $f_{k+a}((y,b))=[|y|-S,|y|]$. Since $x \leq y$, we have $|x| \leq |y|$ and so $f_{k+a}((x,b)) \cap f_{k+a}((y,c))= [|x|,|y|]$. Hence, $\{(x,b),(y,c)\} \in E(I_{k+a})$.
\end{proof}

\begin{claim}\label{Type 2 interval graphs and non-edge}
    For any Type 2 non-edge $\{(x,b),(y,c)\}$ of $\TT$,
    there exists some Type 2 unit interval graph $I_{k+a}$ such that $\{(x,b),(y,c)\} \notin E(I_{k+a})$.
\end{claim}
\begin{proof}
    Since $\{(x,b),(y,c)\}$ is a non-edge of Type 2, without loss of generality, we can assume that $x >y$ and $b < c$. Now $b < c$ implies $b+1 \leq c$. By construction of the interval graph $I_{k+b+1}$, the vertices $(x,b)$ and $(y,c)$ are mapped as follows: $f_{k+b+1}((x,b))=[|x|,S+|x|]$ and $f_{k+b+1}((y,c))=[|y|-S,|y|]$. Since $x > y$ implies $|x| > |y|$, the intervals $[|x|,S+|x|]\cap [|y|-S,|y|] = \emptyset$ and hence, $(x,b)$ and $(y,c)$ are non-adjacent in $I_{k+b+1}$.
\end{proof}

  From \Cref{Type 2 interval graphs are supergraphs} and \Cref{Type 2 interval graphs and non-edge}, one can see that the intersection $\bigcap_{a=1}^{a_s}I_{k+a}$ of the Type 2 unit interval graphs contains each edge of $\TT$, but does not contain any Type 2 non-edge.
  
  From our discussion so far (\Cref{type 1 interval graphs are supergraphs} - \Cref{Type 2 interval graphs and non-edge}), we can conclude that $\TT=\bigcap_{i=1}^{k+a_s}I_i$. Hence, by \Cref{boxicity def 2}, $cub(\TT) \leq k+a_s$. 
\end{proof}

Note that $TCC(a_1)$ is a complete graph, and hence $cub(TCC(a_1))=0$. Therefore, using \Cref{inductive process}, the following holds: $cub(\TT) \leq a_2+a_3+\dots + a_{s}$. Moreover, since the ordering of the complete graphs does not affect their Cartesian product, $\M \cong TCC(m_d,m_{d-1},\dots,m_2,m_1)$. Then putting $s=d$ and $a_1=m_d,a_2=m_{d-1},\dots,a_{d-1}=m_2,a_d=m_1$ in the last mentioned inequality, we have $cub(TCC(m_1,\dots,m_d)) \leq m_{d-1}+m_{d-2}+\dots + m_1$. Hence, \Cref{upper bound} is proved.\hfill $\lhd$

\section{Proof of \Cref{lower bound}}


\textit{Proof plan:} We identify several vertex-disjoint induced subgraphs in the graph $\M$, each of which is isomorphic to $\TTCH{d}$ or $\TTCH{j}$, where $2 \leq j \leq (d-1)$. Then, we show that the vertex sets of all these induced subgraphs together induce the join of these induced subgraphs in $\M$. Since the boxicity of the join is the sum of the boxicity of its constituent parts (\Cref{general join}), a general expression for the lower bound for $box(\M)$ is derived (\Cref{general formula for lower bound}). We first identify $m_1$ induced subgraphs of $\M$, each of which is isomorphic to $\TTCH{d}$. Then, for $1\leq l \leq d-2$, we identify $(m_{l+1}-m_{l})$ induced subgraphs of $\M$, each of which is isomorphic to $\TTCH{d-l}$. After that, \Cref{lower bound} is an immediate consequence of 
\Cref{induced subgraph is a join}, \Cref{general formula for lower bound} and \Cref{lower bound for boxicity of hypercube}.\hfill $\square$

\vspace{0.4cm}

Now we go into the details of the proof.

\vspace{0.4cm}

\noindent \textit{Proof:} In the proof, we call $x$ a \textit{uniform} $s$-tuple if all $x_i$'s are equal; otherwise, we call $x$ a \textit{non-uniform} $s$-tuple. To simplify our discussion further, we first define a new graph.

~~~

\noindent\textbf{Defining $k$-lifted $H_s$: }For two integers $k\geq 1$ and $s\geq 2$, let $\Gamma_s(k)$ be a graph with the vertex set $V(\Gamma_s(k))=\{x=(x_1,x_2,\dots,x_s) \ : \ x_i \in \Z \text{ and } k-1 \leq x_i \leq k, \forall \ i \in [s], $ and $x$ is a non-uniform $s$-tuple$\}$ and two distinct vertices $x$ and $y$ are adjacent if and only if $x < y$. We show that $\Gamma_s(k)$ is isomorphic to $\TTCH{s}$. This motivates us to call the graph $\Gamma_s(k)$ as \textit{$k$-lifted truncated transitive closure of hypercube $H_s$}. But for simplicity, we abbreviate the name to  
\textit{$k$-lifted $H_s$}.

\begin{lemma}\label{iso to hypercube}
     $\Gamma_s(k) \cong \TTCH{s}$, where $k\geq 1$ and $s\geq 2$ are integers.
\end{lemma}
\begin{proof}
    Note that for any vertex $x \in V(\Gamma_s(k))$, for all $i \in [s]$, $x_i$ is either $k-1$ or $k$. We define a mapping $f: V(\Gamma_s(k)) \xrightarrow{} V(\TTCH{s})$ as follows: for a vertex $x \in V(\Gamma_s(k))$ and for all $i \in [s]$, if $x_i=k-1$, then $(f(x))_i=0$ and if $x_i = k$, then $(f(x))_i=1$. It is easy to see that $f$ is a bijection between $V(\Gamma_s(k))$ and $V(\TTCH{s})$.

    To show that $f$ is an isomorphism, take $x,y \in V(\Gamma_s(k))$ such that $\{x,y\} \in E(\Gamma_s(k))$. Without loss of generality, we can assume that $x < y$. Hence, $k-1 \leq x_i \leq y_i \leq k,$ for all $i \in [s]$. Moreover, there exists at least one $t$ in $[s]$, $x_t < y_t $. This implies $0 \leq (f(x))_i \leq (f(y))_i \leq 1$ for all $i \in [s]$, but $(f(x))_t < (f(y))_t$. Therefore, $f(x) < f(y)$ and $\{f(x), f(y) \} \in E(\TTCH{s})$. Similarly, we can show that if $\{f(x),f(y)\}$ is an edge in $\TTCH{s}$, then $\{x,y\}$ is an edge in $\Gamma_{s}(k)$. Hence $\Gamma_s(k) \cong \TTCH{s}$.
\end{proof}


\vspace{0.2cm}

\noindent {\bf Identifying Induced `$k$-lifted $H_d$' Subgraphs in $\mathbf{\M}$:} For $1 \leq j \leq m_1$, we define the subgraph $\Y_d(j)$ as the subgraph of $\M$ induced on the set $V(\Y_d(j))=\{ x=(x_1,x_2,\dots,x_d) \ : \ x_i \in \Z, \ j-1 \leq x_i \leq j \text{ and each $x$ is a non-uniform $d$-tuple}\}$. The reader may convince themselves that each $V(\Y_d(j))$ is indeed a subset of $V(\M)$. By \Cref{iso to hypercube}, each $\Y_d(j)$, $1\leq j \leq m_1$ is isomorphic to $\TTCH{d}$. Hence, for each $1 \leq j \leq m_1$, we have identified an induced subgraph in $\M$ that is isomorphic to a $j$-lifted $H_d$. Moreover, these induced subgraphs are vertex disjoint, meaning $V(\Y_d(s)) \cap V(\Y_d(t)) = \phi$ for all $1 \leq s < t \leq m_1$. To see this, suppose there exists a vertex $x \in V(\Y_d(s)) \cap V(\Y_d(t))$. Then, we must have $x_i = s-1$ or $s$ (since $x \in V(\Y_d(s))$), and $x_i = t-1$ or $t$ (since $x \in V(\Y_d(t))$), for all $i \in [d]$. Since $s < t$, the only possibility is $s=t-1$ and $x_i=s$ for all $i \in [d]$. However, this implies that $x$ is a uniform $d$-tuple, which contradicts our construction of $V(\Y_d(j))$ as being composed solely of non-uniform $d$-tuples. Hence, the subgraphs $\Y_d(j)$ are indeed vertex-disjoint.

~~~

\noindent {\bf Identifying Induced `$k$-lifted  $H_s$' Subgraphs for $\mathbf{s <d}$ in $\mathbf{\M}$:} Note that, for $k > m_1$, it is not possible to identify an induced subgraph isomorphic to $k$-lifted $H_d$ in $\M$. So, we aim to identify an induced subgraph isomorphic to $k$-lifted $H_s$ for $s < d$. For each $1 \leq  l \leq d-2$, if $m_l < m_{l +1}$, we define the subgraph $\Y_{d-l}(m_l + j)$ for $1 \leq j \leq m_{l+1}-m_l$ as the subgraph of $\M$ induced on the set $\{x=(m_1,m_2,\dots,m_l,x_{l+1},\dots,x_d) \ : \ x_i \in \Z, \ m_l + j -1 \leq x_i \leq m_l +j, \ \forall \ l+1 \leq i \leq d$ and each $(x_{l+1},x_{l+2},\dots,x_d)$ is a non-uniform $(d-l)$-tuple$\}$. The reader may convince themselves that each $V(\Y_{d-l}(m_l+j))$ is indeed a subset of $V(\M)$. By \Cref{iso to hypercube}, one can show that each induced subgraph $\Y_{d-l}(m_l + j)$, $1 \leq j \leq m_{l+1}-m_l$, is isomorphic to $\TTCH{d-l}$ (since the first $l$ components of each vertex are fixed). This is true for all $1 \leq l \leq d-2$. Thus, for each $1 \leq l \leq d-2$ and, for each $1\leq j \leq m_{l+1}-m_l$, we identify an induced subgraph in $\M$ isomorphic to a $(m_l+j)$-lifted $H_{d-l}$. 

~~~

\noindent {\bf Vertex-disjointness of the Identified Induced `$k$-lifted $H_s$' Subgraphs:} Here we show that the vertex sets of identified subgraphs are disjoint, i.e., $V(\Y_d(m_l+j)) \cap V(\Y_{d-l-b}(m_{l+b}+k)) = \emptyset$, for $1 \leq l \leq l+b \leq d-2$, $1 \leq j \leq m_{l+1}-m_l$ and $1 \leq k \leq m_{l+b+1}-m_{l+b}$. Since we have already shown that the induced subgraphs $\Y_d(j)$, $1\leq j \leq m_1$, are pairwise vertex disjoint, we can similarly show the case when $b=0$, i.e, the induced subgraphs $\Y_{d-l}(m_{l}+k)$, $1 \leq l \leq d-2$ and $1 \leq k \leq m_{l+1}-m_{l}$ are pairwise vertex disjoint. So, it is enough to assume that $b \geq 1$. Suppose there exists a vertex $x \in V(\Y_{d-l}(m_l + j)) \cap V(\Y_{d-l-b}(m_{l+b}+k))$, where $b \geq 1$ and, $j$ and $k$ are mentioned as before. Then, by construction, we must have $x_i=m_l +j -1$ or $m_l +j$ for all $l+1 \leq i \leq d$ and, $x_i=m_{l+b}+k-1$ or $m_{l+b}+k$, for all $l+b+1 \leq i \leq d$. As $b \geq 1$, the following holds: $m_{l+b} \geq m_{l+1} > m_l$. So, the only possibility is $b=1$ and $m_l + j = m_{l+1} + k -1$. The last equation holds only when $j = m_{l+1}-m_l$ and $k=1$ (Note that both these values are valid.). In that case, $x_i=m_{l+1}$ for all $l+1 \leq i \leq d$. This means that $x$ is a uniform $(d-l)$-tuple of $V(\Y_{d-l}(m_{l+1}-m_l))$, which contradicts our construction.

~~~

\noindent {\bf Collection of all the Identified Induced `$k$-lifted $H_s$' Subgraphs:} We now show that the vertex sets of the induced subgraphs $\Y_d(j)$, for all $1\leq j \leq m_1$ and $\Y_{d-l}(m_l+j)$, for all $1 \leq l \leq d-2$ and $1 \leq j \leq m_{l+1}-m_l$, together induce the join of these subgraphs in $\M$. This result is stated and formally proved in the following lemma.

\begin{lemma}\label{induced subgraph is a join}
     Let $1 \leq l \leq l+b \leq d-2$. Then, the vertex set of $\Y_{d-l}(m_l+j)$ and the vertex set of $\Y_{d-l-b}(m_{l+b}+k)$ together induce a join in $\M$, for each $1 \leq j \leq m_{l+1}-m_l$ and for each $1\leq k \leq m_{l+b+1}-m_{l+b}$ (Note that $b$ can be equal to $0$. But, when $b=0$, we assume that $j \neq k$.).

     
\end{lemma}
\begin{proof}
    
  To prove the lemma,  it is sufficient to prove the following: (1) $V(\Y_{d-l}(m_l+j)) \cap V(\Y_{d-l-b}(m_{l+b}+k)) = \emptyset$, (2) If $x \in V(\Y_{d-l}(m_l+j))$ and $y \in V(\Y_{d-l-b}(m_{l+b}+k))$, then $x \leq y$.

    Note that we have already proved (1) in our previous discussion. So, we prove (2) now.

    \textit{Case 1:} $b \neq 0$, i.e., $b \geq 1$.

    Here $b\geq 1$ implies $m_{l+b} \geq m_{l+1}$. Recall that, $x=(m_1,m_2,\dots,m_l,x_{l+1},x_{l+2},\dots,x_d)$, such that $m_l + j -1 \leq x_i \leq m_l +j$, for all $l+1 \leq i \leq d$.   
    Also, $y=(m_1,m_2,\dots,m_l,m_{l+1},\dots,m_{l+b},y_{l+b+1},\dots,y_d)$, such that $m_{l+b} + k -1 \leq y_i \leq m_{l+b} +k$, for all $l+b+1 \leq i \leq d$.
    Since $j \leq m_{l+1}- m_l$ and $k \geq 1$, we have $x \leq y$ due to the following relation:
    \begin{equation}\label{eq 2}
        x_s \leq m_l +j \leq m_{l+1} \leq m_{l+b} \leq m_{l+b} + k -1 \leq y_s, \text{ for all } s.
    \end{equation}

    \textit{Case 2:} $b=0$

    Since $b=0$, without loss of generality, we can assume that $j<k$. Therefore, \Cref{eq 2} holds in this case also, as $j < k$ implies $j \leq k-1$. Hence, $x \leq y$.\end{proof}

\vspace{0.2cm}

\noindent \textbf{Deriving a General Expression for a Lower Bound of Boxicity of $\mathbf{\M}$}: From our discussion till now, we can say that if $m_1 < m_2 < \dots < m_d$, then in $\M$ there are $m_1$ induced subgraphs which are isomorphic to a lifted $H_d$ and for $1\leq l \leq d-2$, there are $(m_{l+1}-m_{l})$ induced subgraphs which are isomorphic to a lifted $H_{d-l}$. Moreover, due to \Cref{induced subgraph is a join}, the join of all these induced subgraphs is an induced subgraph in $\M$. That means, $\Y = \big( \bigvee_{j=1}^{m_1} \Y_d(j)\big) \vee \big( \bigvee_{l=1}^{d-2}\bigvee_{j=1}^{m_{l+1}-m_l} \Y_{d-l}(m_l+j)\big)$ is an induced subgraph of the graph $\M$. From \Cref{general join} and \Cref{subgraph}, we have $box(M)\geq box(\Y)= m_1 \cdot box(\TTCH{d}) + \sum_{l=1}^{d-2}(m_{l+1}-m_l) \cdot box(\TTCH{d-l})$.  This leads to the following lemma.

\begin{lemma}\label{general formula for lower bound}
    Let $m_1 \leq m_2 \leq \dots \leq m_d$ be $d$ natural numbers, where $d \geq 2$. Let $f(d)$ be the boxicity of $\TCH{d}$. Then, $box(\M) \geq m_1 f(d) + (m_2-m_1) f(d-1) + (m_3-m_2) f(d-2) + \dots + (m_{d-1}-m_{d-2}) f(2)$.
\end{lemma}




\vspace{0.2cm}

\noindent \textbf{Best Known Lower Bound  for $\mathbf{f(s)}$}: The best known lower bound  of $f(s)= box(\TCH{s})$, where $s \geq 2$, follows from the literature (see our discussion in \Cref{intro: box of divisor graph}). But we prove it here explicitly.


\begin{lemma}\label{lower bound for boxicity of hypercube} $box(TC(H_s)) \geq \ceil{\frac{s}{2}}$, where $s \geq 2$.
\end{lemma}


\begin{proof}
    For $s=2$, it is easy to note that the $box(TC(H_2))=1=\ceil{\frac{2}{2}}$. For $s>2$, we define two subsets $A$ and $B$ of $V(TC(H_s))$ as follows:
    \begin{center}
        $A =\{(1,0,\dots,0), (0,1,0,\dots,0),\dots$, $(0,0,\dots,0,1)\}$,\\
        $B =\{(0,1,1,\dots,1),$$(1,0,1,\dots,1),\dots,(1,1,\dots,1,0)\}$,
    \end{center} where the elements of $A$ and $B$ are $s$-tuples. 
    
    For any two distinct vertices $a, a' \in A$, neither $a < a'$ nor $a > a'$ and hence $\{a, a'\}$ $\notin$ $E(TC(H_s))$. So, $A$ is an independent set in $TC(H_s)$. Similarly, $B$ is also an independent set in $TC(H_s)$. Hence, the subgraph $TC(H_s)[A \cup B]$ of $TC(H_s)$, induced on $A\cup B$, is a bipartite graph with $|A| = |B| = s$. 

    Define a bijection $f : A \to B$ such that for each $a \in A$,  
\[
f(a) = b \text{ where } b_i = 
\begin{cases} 
1 & \text{if } a_i = 0, \\ 
0 & \text{if } a_i = 1,
\end{cases}
\]
for all $1 \leq i \leq s$. For example, $f((1, 0, 0, \dots, 0)) = (0, 1, 1, \dots, 1)$.
    
    Note that for each vertex $a \in A$, neither $a < f(a)$ nor $a > f(a)$. Moreover, for each vertex $a \in A$, we have $a < b$, for all $b \in B\backslash\{f(a)\}$. Therefore, $TC(H_s)[A \cup B]$ is isomorphic to a bipartite graph obtained by removing a perfect matching from a complete bipartite graph $K_{s,s}$, where the perfect matching is $\{\{a, f(a)\} \ : \ a \in A\}$. From \cite{chandran2009cubicity}, the boxicity of this bipartite graph is $\ceil{\frac{s}{2}}$. Hence, by \Cref{subgraph}, $box(TC(H_s) \geq \ceil{\frac{s}{2}}$.
 \end{proof}

\vspace{0.2cm}

\noindent \textbf{Deriving the Final Lower Bound Result}: As a corollary of \Cref{general formula for lower bound} and \Cref{lower bound for boxicity of hypercube}, we establish the lower bound of boxicity of $\M$ as stated in \Cref{lower bound}.
From \Cref{general formula for lower bound} and \Cref{lower bound for boxicity of hypercube}, we have the following inequality:
\begin{equation}\label{eq3}
   box(M) \geq m_1 \ceil{d/2} + (m_2-m_1) \ceil{(d-1)/2} + \dots + (m_{d-1}-m_{d-2}) \ceil{2/2}. 
\end{equation}
If $d=2$, \Cref{eq3} implies $box(M) \geq m_1 \ceil{2/2}$. We divide the case $d\geq 3$ into two subcases when $d$ is even and $d$ is odd. Note that when $d$ is even, we have the following:
\begin{center}
    $\ceil{(d-i)/2}=\ceil{(d-i-1))/2}=(d-i)/2$, where $0\leq i \leq (d-2)$ and $i$ is even.
\end{center}
When $d$ is odd, we have the following:
\begin{center}
   $\ceil{d/2} = (d+1)/2$, \\ $\ceil{(d-i)/2}=\ceil{(d-i-1))/2}=(d-i)/2$, where $1 \leq i \leq (d-2)$ and $i$ is odd.
\end{center}

A simple calculation after these substitutions  in \Cref{eq3} gives the required result.\hfill $\lhd$

\section{Open Problems}

If the lower bound for the boxicity of $TC(H_d)$ mentioned in \Cref{lower bound for boxicity of hypercube} is improved, it will consequently improve the lower bound for the boxicity of $\M$ in view of \Cref{general formula for lower bound}. This naturally leads to our first open problem. Furthermore,  we believe that the additional open problems are natural extensions of our work.
\begin{enumerate}
    \item What is the optimal value of $f(d)$ (mentioned in \Cref{general formula for lower bound})? Moreover, what is the exact value of the boxicity of the divisor graph $D(n)$?
    
    \item What is the boxicity of divisor graphs defined over more general vertex sets? In particular, the sets such as \( (n\alpha, n] \), \( a[n] + b = \{ak + b : k \in [n]\} \), and \( (\alpha n, n] \), for \( 0 < \alpha < 1 \), might be of particular interest, since Lewis et al. \cite{lewis2021order} previously investigated the poset dimension of the corresponding divisibility posets.
    
    \item What is the boxicity of power graphs of arbitrary groups?
\end{enumerate}

\bibliographystyle{alpha}
\bibliography{references}

\end{document}